# Markowitz-based cardinality constrained portfolio selection using Asexual Reproduction Optimization (ARO)

**Taha Mansouri,[a], Mohammad Reza Sadeghi Moghadam,[1b],**

[a] Postdoctoral researcher, School of science, engineering, and environment, university of Salford, Email:t.mansouri@salford.ac.uk, Tel:+447411781036

[b] Department of Production and Operation Management, Faculty of Management, University of Tehran, Tehran, Iran, Email:rezasadeghi@ut.ac.ir, Tel:+989126396643

---

Corresponding author[1]

**Abstract:**

   The Markowitz-based portfolio selection turns to an NP-hard problem when considering cardinality constraints. In this case, existing exact solutions like quadratic programming may not be efficient to solve the problem. Many researchers, therefore, used heuristic and metaheuristic approaches in order to deal with the problem. This work presents Asexual Reproduction Optimization (ARO), a model free metaheuristic algorithm inspired by the asexual reproduction, in order to solve the portfolio optimization problem including cardinality constraint to ensure the investment in a given number of different assets and bounding constraint to limit the proportions of fund invested in each asset. This is the first time that this relatively new metaheuristic is in the field of portfolio optimization, and we show that ARO results in better quality solutions in comparison with some of the well-known metaheuristics stated in the literature. To validate our proposed algorithm, we measured the deviation of obtained results from the standard efficient frontier. We report our computational results on a set of publicly available benchmark test problems relating to five main market indices containing 31, 85, 89, 98, and 225 assets. These results are used in order to test the efficiency of our proposed method in comparison to other existing metaheuristic solutions. The experimental results indicate that ARO outperforms Genetic Algorithm(GA), Tabu Search (TS), Simulated Annealing (SA), and Particle Swarm Optimization (PSO) in most of test problems. In terms of the obtained error, by using ARO, the average error of the aforementioned test problems is reduced by approximately 20 percent of the minimum average error calculated for the above-mentioned algorithms.

**Key Words:** Portfolio optimization, Cardinality constraints, Markowitz mean-variance model, Asexual reproduction optimization, Efficient frontier.

1. **Introduction**:

Portfolio optimization, which is the problem of allocating the initial amount of capital among a given number of assets or securities, has attracted a lot of attention in the field of quantitative finance (Moral-Escudero, Ruiz-Torrubiano, & Suarez, 2006). In order to help investors in optimally forming their portfolio of assets, Markowitz (1952, 1959) has proposed a quantitative framework. Markowitz mean-variance model of portfolio selection which caused the development of *Modern Portfolio Theory (MPT)*, formulates the problem as a multi-objective optimization problem with two conflicting objectives: maximizing the expected return and minimizing the risk

(measured by variance) of a portfolio. Considering these two competing criteria simultaneously, there is no single optimal portfolio, but a set of portfolios forming the *efficient frontier (EF)*. In other words, efficient frontier, also called *Pareto-optimal front*, is the collection of portfolios which result in minimum risk for a given level of return or equivalently, maximum return for a given level of risk.

The standard Markowitz model assumes a perfect market without any transaction costs and taxes, where short selling is forbidden, and assets are tradable in any non-negative fractions. This basic model belongs to the category of Quadratic Programming (QP) problems (Fernández & Gómez, 2007); thus, the efficient frontier could be found using standard QP solvers which are easily available and can guarantee to find the optimal solution, and be modified to include linear constraints (Chang, Meade, Beasley, & Sharaiha, 2000). However, in the absence of these unrealistic assumptions, or presence of some non linear real-world constraints, QP is not necessarily feasible for finding efficient portfolios any more.

Many researchers have tried to extend the Markowitz's model (Markowitz, 1952) in order to capture more realistic market conditions by introducing some additional constraints. These include, cardinality constraint which limits the number of assets held in the portfolio, bounding (quantity or floor-ceiling) constraint, also known as buy-in thresholds, imposing lower and/or upper bounds on funds invested in each asset, pre-assignment constraint reflecting the investor's preferences by requiring some specific assets to be held in the portfolio, round lot (minimum lots) constraint which forces the amount invested in each asset to be a multiple of *minimum transaction lot*, class constraint which limits the total weight assigned to a class (assets with common characteristics), as well as turnover and trading constraints that impose upper and lower bounds respectively on the variation of the assets weight from one period to another, which are particularly useful in multi period investments (Ponsich & Antonio, 2012; Lwin et al., 2014; Crama & Schyns, 2003; Tollo & Roli, 2008).

According to Metaxiotis & Liagkouras (2012), cardinality and bounding constraints occupy the main focus of the researchers. In practice, many investors prefer to hold a certain number of assets in their portfolio so as to facilitate its management, decrease transaction costs, and assure a minimum degree of diversification. Moreover, they prefer to avoid holding very small and large proportion of assets in order to reduce administrative costs and risk, respectively (Anagnostopoulos & Mamanis, 2011; Lwin et al., 2014).

In this paper, we tackle the problem with regard to the extended Markowitz mean-variance model which includes cardinality and bounding constraints. In this case, we refer to the EF as cardinality constrained efficient frontier (CCEF). By introducing the cardinality constraint into the classic quadratic programming model, this problem turns to a mixed integer quadratic programming one which is NP-hard (Bienstock, 1996; Moral-Escuderoet al., 2006; Shaw, Liu, & Kopman, 2008). In this case, exact optimization methods are not efficient for large problem sizes (Kalayci et al., 2017). Many researchers, therefore, take advantage of heuristic and metaheuristic approaches in order to deal with the problem (Maringer, 2006; Tollo & Roli, 2008). Although these approaches do not guarantee to find the optimal solution, they are efficient for finding near-optimal solutions.

Asexual Reproduction Optimization (ARO), proposed relatively recently in Farasat et al. (2010) and Mansouri et al. (2011), is an evolutionary individual based metaheuristic algorithm inspired by *budding* mechanism of asexual reproduction and has been used in very few studies (e.g. Khanteymoori et al., 2011; Kazemi et al., 2012; Noormohammadi Asl et al., 2014; Ahmadian & Khanteymoori, 2015; Yazdanparast et al., 2015). None of these studies deal with the portfolio selection problem (PSP).

The ARO has advantages that make it completely different from other metaphors. First, it is an individual-based algorithm. Thus, unlike population-based algorithms that require a large amount of computational resources to convert, ARO consumes much less. Hence, it converges faster. The second case is mathematical convergence, so it has good exploration and exploitation rates. Third, ARO does not require parameter settings, so you are unlikely to have trouble setting parameters that are a common meta-cognitive problem such as genetic algorithms (GA), annealing simulation (SA), taboo search (TS). And Particle Particle Optimization (PSO). In addition, the ARO does not use any selective mechanism such as a roulette wheel. Inappropriate selection of selection mechanisms may lead to problems such as premature convergence due to excessive selection pressure. Fourth providers in many benchmark issues have shown the computational power of this algorithm. Fifth, ARO is a free model algorithm that can be applied to various types of optimization. (Mansouri et al., 2011) Finally, the ARO in this paper presents better results than the algorithms used in other papers.

For these reasons, we take advantage of ARO to tackle the Markowitz-based cardinality constrained portfolio selection problem. The main contribution of this study is to solve this problem more efficiently using a new approach. We apply a method which uses ARO to confront the portfolio selection problem. Our proposed method results in better quality solutions compare to some of the well-known metaheuristics which have been used in this field.

Computational results are reported for five analyses of weekly price data with regard to the following indices for the time period between March, 1992 to September, 1997: Hang Seng 31 in Hong Kong, DAX 100 in Germany, FTSE 100 in the UK, S&P 100 in the USA and Nikkei 225 in Japan.

The rest of this paper is organized as follows. A literature survey of exact, heuristic, and hybrid approaches to the problem is presented in Section 2. Section 3 describes the generic mean-variance portfolio selection problem, followed by the specific model in the presence of cardinality and bounding constraints. Section 4 introduces the proposed ARO algorithm along with its application to the problem under consideration. Computational experiments and results are discussed in Section 5. Conclusion and future work are presented in Section 6.

## 2. Literature review

According to Woodside-Oriakhi, Lucas, & Beasley (2011), researchers dealt with the cardinality constrained portfolio optimization using either exact or heuristic approaches. In this paper, we consider a third category to structure our literature survey: hybrid methods which combine an exact method with a heuristic one.

*2.1. Exact approaches*

As stated earlier, when considering the cardinality constraint into the model, exact methods may not be efficient to solve portfolio optimization problem for large problem samples. However, some researchers tried to deal with the problem using a relaxed version of cardinality constraint which imposes an upper bound on the number of assets present in the portfolio. This approach, in which Eq. (14) is an inequality rather than equality, has a significantly less computational complexity (Woodside-Oriakhi et al., 2011). Moreover, the results show that researchers were able to handle this version of problem for limited problem sizes (Lwin et al., 2014).

Table 1 summarizes the exact approaches used in the literature to solve the PSP.

**Table 1**
Exact approaches for portfolio selection problem.

| Authors | Constraints | Method | Datasets | risk measure |
|---|---|---|---|---|
| Bienstock (1996) | upper bound on the cardinality | branch-and-cut | 3897 assets | variance |
| Lee & Mitchell (1997) | upper bound on the cardinality | interior point nonlinear branch-and-bound solver | 5 datasets involving up to 150 assets | variance |
| Young (1998) | transaction costs | minimax | No computational results | minimum return |
| Li et al. (2006) | upper bound on the cardinality, round lots | convergent Lagrangian and contour-domain cut | 30 assets from from the Hong Kong stock market | Variance |

| Shaw et al. (2008) | upper bound on the cardinality | Lagrangian relaxation | 8 datasets involving up to 500 assets | Variance |
|---|---|---|---|---|
| Vielma et al. (2008) | upper bound on the cardinality | LP based branch-and-bound | Problem sets involving up to 200 assets | Variance |
| Bertsimas & Shioda (2009) | upper bound on the cardinality, bounding | tailored branch-and-bound | Problem sets involving up to 500 assets | variance |
| Gulpinar et al. (2010) | cardinality, bounding, round lots | DC functions programming | 98 assets | variance |

## 2.2. *Heuristic approaches*

With regard to the above-mentioned approach, using other risk measures than variance, Mansini & Speranza (1999) considered PSP with minimum transaction lots and showed that in this case, the problem of finding a feasible solution is NP-complete no matter what the risk measure is. In their work, they used Mean Semi-absolute Deviation as a measure of risk and presented three heuristics based on solving the linear programming relaxation to tackle the problem. Kellerer, Mansini, & Speranza (2000) also considered the same risk measure in their paper. They added fixed transaction costs to the previous model and employed two of the three Mixed Integer Linear Programming based (MILP-based) heuristics which were used in the previous study. In a more recent work, Chang, Yang, & Chang (2009) considered different risk criteria other than variance; semi-variance, mean absolute deviation (MAD) and variance with skewness. They employed Genetic Algorithm (GA), and showed its efficiency for solving these problems in different risk measures.

Heuristic attempts to solve the portfolio selection problem are summarized in Table 2.

**Table 2**
Heuristic approaches for portfolio selection problem.

| Authors | Constraints | Method | Datasets | risk measure |
|---|---|---|---|---|
| Mansini & Speranza (1999) | minimum transaction lots | Basic, reduced cost, and iterated MILP-based heuristics | 2 datasets involving 244 and 277 assets | Mean Semi-absolute Deviation |
| Kellerer et al. (2000) | fixed transaction costs, round lots | reduced cost and iterated MILP-based heuristics | 1 dataset involving 244 assets | Mean Semi-absolute Deviation |
| Chang et al. (2000) | cardinality, bounding | GA, SA, TS | Hang Seng, DAX, FTSE, S&P100, Nikkei | variance |
| Jobst et al. (2001) | cardinality, bounding, round lots | integer restart, reoptimisation | test problems provided by Chang et al. (2000) | variance |
| Schaerf (2002) | cardinality, bounding | hill climbing, SA, TS | test problems provided by Chang et al. (2000) | variance |
| Crama & Schyns (2003) | Cardinality, turnover, trading | SA | 1 dataset involving 151 US stocks | variance |

| Author | Constraints | Method | Data | Risk measure |
|---|---|---|---|---|
| Derigs & Nickel (2003) | Considering any constraint is possible in this model | a metaheuristic approach based on SA | a case study containing 202 stocks related to DAX30 & STOXX200 | variance |
| Maringer & Kellerer (2003) | upper bound on the cardinality | a combination of SA and ES | DAX30, FTSE100 | variance |
| Ehrgott et al. (2004) | cardinality, bounding | customized local search, SA, TS, GA | four data sets involving up to 1416 assets | volatility, S&P star ranking |
| Fernández & Gómez (2007) | cardinality, bounding | HNN | test problems provided by Chang et al. (2000) | variance |
| Chiam et al. (2008) | cardinality, bounding | an approach based on a multiobjective EA | test problems provided by Chang et al. (2000) | variance |
| Chang et al. (2009) | cardinality, bounding | GA | HANG SENG, FTSE, S&P100 | semi-variance, MAD, variance with skewness |

| Pai & Michel (2009) | Cardinality, bounding, class | an ES-based solution and k-means clustering | BSE 200, Nikkei 225 | variance |
| --- | --- | --- | --- | --- |
| Soleimani et al. (2009) | cardinality, round lots, market (sector) capitalization | GA | two randomly generated datasets with 500 and 2000 assets. | variance |
| Cura (2009) | cardinality, bounding | PSO | test problems provided by Chang et al. (2000) | variance |
| Anagnostopoulos & Mamanis (2010) | class, bounding | NSGA-II, PESA, SPEA2 | two randomly generated data sets containing 200 and 300 assets | variance |
| Anagnostopoulos & Mamanis (2011) | cardinality, bounding | NPGA2, NSGA-II, PESA, SPEA2, e-MOEA | seven test problems involving up to 2196 assets | variance |
| Woodside-Oriakhi et al. (2011) | cardinality, bounding | GA, SA, TS | seven test problems involving up to 1318 assets | variance |

| Deng et al., 2012) | cardinality, bounding | PSO | test problems provided by Chang et al. (2000) | variance |
| --- | --- | --- | --- | --- |
| Lwin & Qu (2013) | cardinality, bounding | PBILDE | test problems provided by Chang et al. (2000) | variance |
| Lwin et al. (2014) | cardinality, bounding, pre-assignment, round lots | MODEwAwL | seven test problems involving up to 1318 assets | variance |
| Ni et al. (2017) | cardinality, bounding | DRTCPSO-AD, DRTCPSO-D, DRTCPSO-LIAD, DRTCPSO-LID | test problems provided by Chang et al. (2000) | variance |
| Kalayci et al. (2017) | cardinality, bounding | ABC | Hang Seng, DAX, FTSE, S&P100, Nikkei, XU030, XU100 | variance |
| Meghwani & Thakur (2017) | cardinality, bounding, pre-assignment | NSGA-II, MOEA/D, GWAS-FGA | FF38, FF48 | VaR, CVaR |

| | , self-financing | | | |
|---|---|---|---|---|
| Liagkouras, K., & Metaxiotis, K. (2018) | Liquidity Upper and Lower limit of the cardinality constraint lower and upper bound constraints, transaction costs | proposed new Multi-period Fuzzy Portfolio Optimization Algorithm (MFPOA) | 92 assets from FTSE-100 index in London | possibilistic variance of the portfolio return at each period. |
| Babazadeh, H., & Esfahanipour, A. (2019) | Cardinality Constraint Budget Constraint at each period | new design of NSGA-II | 20 selected assets from S&P100 | VaR |
| Ling, A., Sun, J., & Wang, M. (2019) | stochastic constraint Budget Constraint | A computationally tractable approximation approach based on second-order cone optimization | 10Ind, 12Ind, iShare and DJIA, as in Rujeerapaiboon et al. (2016) , with all data from the French data library, Yahoo finance and | lower partial moment (LPM). |

| | | | the Wind database | |
|---|---|---|---|---|
| Ahmadi-Javid, A., & Fallah-Tafti, M. (2019) | cardinality, bounding | primal-dual interior-point algorithm EVaR- PD algorithm | 50 assets 100 assets 200 assets 750 assets 1000 assets | entropic value-at-risk (EVaR) |
| Gupta, P., Mehlawat, M. K., Yadav, S., & Kumar, A. (2019) | Budget Constraint, Cardinality No short selling constraint Contingent constraint | polynomial goal programming approach | a sample of twenty assets (n = 20) from National Stock Exchange (NSE), India | Variance entropy |
| Mohammadi, S., & Nazemi, A. (2020) | Budget Constraint | neural network | 10 securities | VaR |

### 2.3. Hybrid approaches

More recently in the literature, researchers tried to tackle the portfolio optimization problem by implementing hybrid strategies that take advantage of both exact and heuristic approaches. Table 3 summarizes the Hybrid methods proposed in the literature.

**Table 3**
Hybrid approaches for portfolio selection problem.

| Research | Constraints | Method | Datasets | risk measure |
|---|---|---|---|---|
| (Moral-Escudero et al. (2006) | cardinality, bounding | combination of GA and QP | test problems provided by Chang et al. (2000) | variance |
| Streichert & Tanaka-yamawaki (2006) | upper bound on the cardinality, maximum limits with one exception | combination of a MOEA and a QP local search | Hang Seng 31, DAX 100 | variance |
| Branke et al. (2009) | Cardinality, bounding, 5-10-40 rule from German investment law | envelope-based MOEA | Hang Seng 31, S&P 100, Nikkei 225, and a benchmark with 500 assets | variance |
| Ruiz-torrubiano & Suárez (2010) | upper bound on the cardinality, bounding | Preprocessing (pruning), SA, GAs, EDAs | test problems provided by Chang et al. (2000) | variance |
| Baykasoğlu et al. (2015) | cardinality, bounding | GRASP-QUAD | test problems provided by Chang et al. (2000) | Variance |
| Li, B., Zhu, Y., Sun, Y., Aw, G., & Teo, K. L. (2018) | Bankruptcy Budget Constraint transaction costs Total wealth at each period | GA combined with penalty function | 8 stocks | Variance |
| Zhang, J., & Li, Q. (2019) | Liquidity risk-free asset constrained in each period lower and upper limits No short selling of assets Cardinality | Hybrid DA-GA (Dragonfly Algorithm-Genetic Algorithm) | 10 risky assets | semi-entropy |

*2.4. Comment*

According to the literature reviewed above, since the pioneering work of Harry Markowitz (Markowitz, 1952), the mean-variance model of portfolio optimization has been the main framework for choosing optimal portfolios. Some extensions have been proposed for this model and among them, MVCCPO including cardinality and bounding constraints has attracted the most of researchers' attention. Using metaheuristic algorithms became the main trend for dealing with this model after the research conducted by Chang et al. (2000) and hybrid methods which take advantage of both heuristic and exact solutions, have been used more recently in the literature.

The contribution of our paper is to complement the reviewed literature by proposing a new approach for portfolio selection based on the Markowitz mean-variance-model which results in better quality solutions in comparison with some of the well-known metaheuristics stated in the literature.

**3. Problem Formulation**

Let us start with the basic (unconstrained) Markowitz model. In its multiobjective form, it can be formulated as follows:

$$\text{Minimize} \sum_{i=1}^{N} \sum_{j=1}^{N} w_i w_j \sigma_{ij} \qquad (1)$$

$$\text{Maximize} \sum_{i=1}^{N} w_i \mu_i \qquad (2)$$

$$\text{subject to} \sum_{i=1}^{N} w_i = 1 \qquad (3)$$

$$0 \leq w_i \leq 1, \quad i = 1, \dots, N \qquad (4)$$

Where $N$ is the number of available assets, $\mu_i$ is the expected return of asset $i$, $\sigma_{ij}$ is the covariance between asset $i$ and $j$, and $w_i$ is the decision variable representing the proportion of money invested in asset $i$. Eq. (1) minimizes the risk of the portfolio (measured by variance) while Eq. (2) maximizes the expected return of the portfolio. Eq. (3) defines the budget constraint which forces the investment of all the money in hand, i.e., asset weights must sum up to one. Finally, Eq. (4) states that all weights should be nonnegative.

By solving the above model, a set of efficient portfolios can be found. These Pareto-optimal (non-dominated) solutions form the unconstrained efficient frontier (UEF), i.e., a continuous curve representing the best possible trade off between risk and return.

This bi-objective model can be also represented as a single objective optimization problem; therefore, it could be solved by applying single objective solution techniques. The famous single objective representation of the basic Markowitz model is as follows:

$$\text{Minimize} \sum_{i=1}^{N}\sum_{j=1}^{N} w_i w_j \sigma_{ij} \qquad (5)$$

$$\text{subject to} \sum_{i=1}^{N} w_i \mu_i = R^* \qquad (6)$$

$$\sum_{i=1}^{N} w_i = 1 \qquad (7)$$

$$0 \leq w_i \leq 1, \quad i = 1, \dots, N \qquad (8)$$

This model attempts to minimize risk by considering the expected return as a constraint. Hence, solving the above single objective problem for different levels of expected return results in tracing the unconstrained efficient frontier.

According to Chang et al. (2000), designing a heuristic based on the above formulation is difficult in that it requires the expected return of the portfolio to be exactly $R^*$.

In practice, for tracing the UEF, a popular approach is to introduce a weighting parameter $\lambda$ ($0 \leq \lambda \leq 1$); thus, the objective function could be represented in a Lagrangian relaxation form (Chang et al., 2000; Anagnostopoulos & Mamanis, 2011):

$$\text{Minimize } \lambda \sum_{i=1}^{N}\sum_{j=1}^{N} w_i w_j \sigma_{ij} - (1-\lambda) \sum_{i=1}^{N} w_i \mu_i \qquad (9)$$

$$\text{subject to} \sum_{i=1}^{N} w_i = 1 \qquad (10)$$

$$0 \leq w_i \leq 1, \quad i = 1, \ldots, N \tag{11}$$

By solving this QP problem for various values of λ, the UEF can be traced from the portfolio with maximum return (λ = 0) to the portfolio with minimum level of risk (λ = 1). Chang et al. (2000) showed that when considering the unconstrained problem, by varying λ in Eq. (9), we can obtain exactly the same efficient frontier as we would get by solving Eqs. (5) - (8) for varying values of R*.

In order to find the cardinality constrained efficient frontier (CCEF), many researchers extended the above-mentioned model by adding cardinality and bounding constraints (e.g. Chang et al., 2000; Fernández & Gómez, 2007; Cura, 2009; Woodside-Oriakhi, Lucas, & Beasley, 2011; Deng, Lin, & Lo, 2012; Baykasoğlu, Yunusoglu, & Burcin Özsoydan, 2015):

$$Minimize\ \lambda \sum_{i=1}^{N} \sum_{j=1}^{N} w_i w_j \sigma_{ij} - (1-\lambda) \sum_{i=1}^{N} w_i \mu_i \tag{12}$$

$$subject\ to\ \sum_{i=1}^{N} w_i = 1 \tag{13}$$

$$\sum_{i=1}^{N} z_i = K \tag{14}$$

$$\varepsilon_i z_i \leq w_i \leq \delta_i z_i, \quad i = 1, \ldots, N \tag{15}$$

$$z_i \in \{0,1\}, \quad i = 1, \ldots, N \tag{16}$$

Where $z_i$ is the decision variable indicating the existence of each asset in the portfolio, hence it is equal to 1, if asset $i$ is included in the portfolio and zero otherwise. Eq. (14) defines the cardinality constraint (portfolio consists of exactly $K$ assets) and Eq. (15) defines the bounding constraint which imposes lower and upper limits on the weight of each asset. In this work, we will consider the same MVCCPO model (Eqs. (12) – (16)).

## 4. ARO for portfolio selection problem

In this section, we present our proposed algorithm for solving the cardinality constrained portfolio selection problem. First, we give a brief overview of the general ARO, then the particular implementation of this proposed algorithm that is customized for finding the CCEF will be presented.

*4.1. Asexual Reproduction Optimization*

ARO, which is an individual based evolutionary algorithm modelling the budding mechanism of asexual reproduction, was first described by Farasat et al. (2010) and Mansouri et al. (2011). In ARO, each individual produces a *bud* via a reproduction mechanism; afterward, the bud and its parent compete with respect to their fitness which is obtained from the objective function of the underlying optimization problem. Through competition for limited resources, the fitter one will survive, while the other will be discarded. This reproduction cycle is repeated until the stopping criteria are met.

Consider the following optimization problem:

$$Maximize\ f(X) \quad subject\ to\ X \in S \tag{17}$$

where $X = (x_1, x_2, ..., x_n); x_i \in \mathbb{R}\ (i = 1, 2, ..., N)$ are the decision variables, $f(X)$ is the objective function, and $S$ defines the search space.

The pseudo code of ARO is illustrated in Figure 1.

**FIGURE 1 ABOUT HERE**

*4.2. The proposed ARO for finding the CCEF*

Here, we introduce the customized version of ARO to deal with the cardinality constrained portfolio optimization problem.

### 4.2.1. Notation

Table 4 introduces the notations used to describe the proposed version of ARO to solve the problem.

**Table 4**
Notations for the proposed ARO.

| Symbol | Description |
|---|---|
| $N$ | the whole number of available assets |
| $K$ | the number of assets present in the portfolio |
| $w_i$ | the proportion of capital invested in asset $i$ |
| $\varepsilon_i$ | the lower bound on the proportion invested in asset $i$ |
| $\delta_i$ | the upper bound on the proportion invested in asset $i$ |
| $T_{max}$ | the number of running iterations for ARO |
| $R[x, y]$ | a random integer number in $[x, y]$ |
| $r[0, 1]$ | a random real number in $[x, y]$ |
| $R$ | the set of $i$ whose proportions are fixed at $\delta_i$ |
| $Q$ | the set of $K$ distinct assets in the current solution |
| $g$ | the length of selected substring from parent's *chromosome* |
| $b$ | the number of buds reproduced from the current parent |

### 4.2.2. Solution representation and encoding

In our solution representation, a vector of size *2K* is used to represent a portfolio. This vector consists of two distinct parts, the first part indicates the asset indices present in the portfolio, and the second part determines the proportion of capital to be invested in each asset in the portfolio. So, the first part would be an integer vector of size *K* with its elements belonging to $\{1, 2, ..., N\}$, and the second part consists of *K* real numbers from $[0,1]$ (Fig. 2).

**FIGURE 2 ABOUT HERE**

Where $x_i$ is an integer variable that belongs to $\{1, 2, ..., N\}$ and it represents the index of the $i^{th}$ asset in our portfolio. As mentioned before, we will have K distinct assets in our portfolio, so, $(i = 1, 2, ..., K)$. In the second part of the solution representation, $w_i \in [0,1]$ denotes the value of investment in the asset $x_i$. For instance, if K= 4, which means that we are constrained to have 4 distinct assets in our portfolio, and the whole number of assets in the market is N=100, we should pick 4 distinct integers from $\{1, 2, ..., 100\}$ to represent the asset indices we are going to hold in the portfolio. Assume that we pick asset number 1, 7, 34, and 87, and put the equal weights for them in the portfolio. Our solution representation will be like this:

**FIGURE 3 ABOUT HERE**

It is important to note that each number in $\{1, 2, ..., N\}$ cannot appear more than once in the integer part of the chromosome.

### 4.2.3. *Constraints satisfaction*

To meet the bounding constraint and ensuring that the sum of the proportions invested in assets equals one (Eq. (15), and Eq. (13)), the following approach is applied based on Chang et al. (2000):

Let $Q$ be the set of $K$ distinct assets in the current solution. The lower limits constraint can be satisfied if all weights in the candidate solution are adjusted by setting $w'_i = \varepsilon_i + w_i(1 - \sum_{i \in Q} \varepsilon_i) / \sum_{i \in Q} w_i$. Note here that $\sum_{i \in Q} w'_i = 1$. So, the above formula satisfies both the lower proportion limits and sum to one.

Let $R$ be the set of $i$ whose proportions are fixed at $\delta_i$. In order to satisfy upper limits constraint, an iterative algorithm can be applied as shown in Figure 4.

**FIGURE 4 ABOUT HERE**

It is also noteworthy that there is no need for any mechanism to handle the cardinality constraint (Eq. (14)), since our solution representation presented in Section 4.2.2 requires each solution to contain exactly $K$ distinct assets.

### 4.2.4. Mutation

In order for our proposed ARO to maximize diversity, we present two types of mutation.

#### 4.2.4.1. Mutation of shares

In this type of mutation, a bud is reproduced by altering some *genes* from the integer part of its parent's chromosome. In optimization terms, a new solution is generated by changing some asset indices present in the portfolio while weights remain unchanged.

In order to reproduce the bud, a substring of length $g$ is randomly selected from the integer part of the parent's chromosome. Thereafter, the genes presented in this substring are replaced with $g$ integers which are absent in the remaining string.

To clarify more, suppose that we have $N = 10$, $K = 5$, and the integer part of the parent's chromosome is shown in Figure 5. In order to select a substring, we randomly generate two distinct integers in $[1, K]$, i.e., $r_1 = [1, K]$ and $r_2 = R[r_1, K]$. so, $g = r_2 - r_1 + 1$. Assume that $r_1 = 2$, $r_2 = 3$. Hence, $\{8, 5\}$ should be eliminated from parent's chromosome, and a string of length 2 composing of two distinct elements from $\{1, 3, 5, 6, 8, 9, 10\}$ (e.g., 3, 8) will be substituted. (see Fig. 5).

**FIGURE 5 ABOUT HERE**

#### 4.2.4.2. Mutation of weights

We use two kinds of variation here to expand the search space, stochastic and chaotic. The former is used when we are stuck in local optimum and helps us to exit from it. Buy using The latter, we try to visit

maximum number of points near the solution. To select the kind of variation, we use the following function:

$$f(i,b) = \sin\left(\max\left(1 - \varphi^{\ln(i)}/b, 0\right) * \frac{\pi}{2}\right) \tag{18}$$

Where $i$ shows the number of iteration, $b$ is a variable representing the number of iterations we searched in local- i.e., the number of buds reproduced from the current parent- and $\varphi$ is the golden number which is approximately equal to 1.618. the value of $f(i,b)$ is decreasing in $i$ and increasing in $b$. Eq. (18) tries to produce a probable measure to determine the extent of being stuck in local optimum. We examined the above-mentioned function to produce this measure and found it sui for our purpose.

Then, base on the value of $f(i,b)$, we use the following procedure to select the variation: A random real number is generated in $[0, 1]$, i.e., $r_3 = r[0, 1]$. if $r_3$ is lower than $f(i,b)$, stochastic variation will be selected. Otherwise, we apply the chaotic variation.

In stochastic variation, we select a random substring of length $g$ from the real part of the parent's chromosome, then define $p$ as follows:

$$p = 1 / (1 + \ln(g)) \tag{19}$$

For each gene from the selected substring, let $r_4, r_5 = r[0, 1]$. If $r_4 \leq p$ and $r_5 \leq 0.3$, the value of the gene will be replaced with $p$ multiplied by a random number in $[0, 1]$. If $r_4 \leq p$ and $r_5 > 0.3$, the new value of the gene takes a random real number in $[0, 1]$. Otherwise, the value of the gene remains unchanged.

In chaotic variation, we apply these steps for every gene in the real part of the parent's chromosome: let $r_6 = r[0, 1]$. If $r_6 \leq 0.2$, the value of the gene should be multiplied by $0.2 * f(i,b)$. If $r_6$ belongs to $[0.3, 0.7]$, then let $r_7 = r[0, 1]$, and the value of the gene will be multiplied by $r_7 + 0.2 * f(i,b)$. Otherwise, the value of the gene remains unchanged.

### 4.2.5. Termination

Our proposed ARO terminates after running for a predefined number of iterations, $T_{max}$.

## 5. Computational results

In this section, our proposed ARO is evaluated and compared to other well-known existing heuristics used for tackling the cardinality constrained portfolio optimization based on the standard test problems. The heuristics which are used for comparison are Genetic Algorithm(GA), Simulated Annealing (SA), Tabu Search (TS), and Particle Swarm Optimization (PSO). We use the results reported by Chang et al. (2000) for GA, SA, and TS while PSO results was those reported by Deng et al. (2012). We report the computational results for finding 50 different portfolios on the CCEF for each data set using the values $\varepsilon_i = 0.01, \delta_i = 1$ ($i = 1, 2, ..., N$), and $K = 10$. We set the number of iterations which is the termination condition for ARO to 20000. The proposed algorithm is implemented in MATLAB language.

### 5.1. Datasets

The performance of our proposed algorithm is evaluated on the benchmark data related to five well-known major market indices from the publicly available OR-Library (Beasley, 1990). These test problems were built based on weekly prices between March 1992 and September 1997 for the following market indices: Hang Seng 31 in Hong Kong, DAX 100 in Germany, FTSE 100 in UK, S&P 100 in USA, and Nikkei 225 in Japan. The number of assets, *N*, related to each dataset is 31, 85, 89, 98 and 225, respectively. These data files contain mean return of each stock, covariance between these stocks, and the unconstrained efficient frontier composing of 2000 points (i.e., *standard efficient frontier*) , which are accessible at http://people.brunel.ac.uk/~mastjjb/jeb/orlib/portinfo.html. They were first provided by Chang et al. (2000) and were used in many other studies since then (e.g. Baykasoğlu et al., 2015; Chiam et al., 2008; Cura, 2009; Deng et al., 2012; Fernández & Gómez, 2007; Lwin & Qu, 2013; Moral-Escudero et al., 2006; Ruiz-torrubiano & Suárez, 2010; Schaerf, 2002).

### 5.2. Performance indicator

To evaluate the performance of a heuristic, the quality of results could be measured in terms of the deviation of obtained results from the optimal solution (Woodside-Oriakhi et al., 2011). In the case of finding the cardinality constrained efficient frontier, because of unavailability of the optimal CCEF, the quality of results could be measured according to their deviation from UEF which can be found simply by QP. Thus, we used exactly the same approach previously proposed by Chang et al. (2000) which is the most commonly used approach in the literature. For instance, the following studies used the same approach: Woodside-Oriakhi et al. (2011), Deng et al. (2012), Lwin & Qu (2013), Baykasoğlu et al. (2015).

Let $(s_p, R_p)$ be the standard deviation and return corresponding to a portfolio $p$ found by ARO heuristic. By using linear interpolation, we can find $s_p^*$ which is the standard deviation associated with $R_p$ in standard efficient frontier. Hence, the *standard deviation error* of portfolio $p$ is defined as follows:

$$\text{Standard deviation error } (p) = 100\left|(s_p - s_p^*)/s_p^*\right| \qquad (20)$$

Similarly, let $R_p^*$ be the return associated with $s_p$ using linear interpolation in standard efficient frontier. Then, the *return error* of portfolio $p$ would be:

$$\text{Return error } (p) = 100\left|(R_p - R_p^*)/R_p^*\right| \qquad (21)$$

Furthermore, the minimum between two above-mentioned errors for portfolio $p$ is defined as *percentage error*, and by averaging this for all obtained portfolios, we can define *mean percentage error*.

$$\text{Percentage error } (p) = \min\{100\left|(s_p - s_p^*)/s_p^*\right|, 100\left|(R_p - R_p^*)/R_p^*\right|\} \qquad (22)$$

$$\text{Mean percentage error} = \sum_{p=1}^{50} percentage\ error\ (p)/50. \qquad (23)$$

In other words, after obtaining all portfolios based on the targeted portfolios, we measure the vertical and horizontal distances to the standard UEF for each portfolio and compute the minimum of these two numbers. Then, our performance indicator, *mean percentage error,* is the average of these derived minimums for all obtained portfolios. For more details about this approach, see Chang et al. (2000).

*5.3.* *Experiments*

As mentioned before, the standard efficient frontier is the set of 2000 optimal portfolios on the UEF which is available from OR-Library (Beasley, 1990) for each of the five tested data sets. Figure 6 shows the heuristic frontier which is formed by the 50 portfolios found by ARO as well as the standard efficient frontier for each data set.

**FIGURE 6 ABOUT HERE**

This figure clearly illustrates that our proposed algorithm performs really well when dealing with portfolios requiring lower levels of risk and expected return. However, when finding portfolios which have much higher risk and expected return, the distance between standard efficient frontier and ARO heuristic frontier increases. The reason is that when we want to constitute portfolios with higher expected return, we should choose fewer assets which have higher mean return; in the unconstrained problem it would be possible to choose even one asset with the highest level of mean return, but when dealing with the constrained problem, our algorithm is forced to select exactly 10 assets because of the cardinality constraint. Hence, the portfolios which are close to top right corner of CCEF, have much significant percentage error

Table 5 compares our results with those obtained by Chang et al. (2000) and Deng et al. (2012) based on the mean percentage error, for each data set. The best mean percentage error among them for each problem is written in boldface. These results show that our proposed method outperforms other heuristics in four out of five test problems. Therefore, the superiority of it is clear from the experimental results.

**Table 5**
Comparison of proposed ARO with other heuristics for finding the CCEF.

| Index | N | GA | SA | TS | PSO | ARO |
| --- | --- | --- | --- | --- | --- | --- |

|         |     | chang et al. | chang et al. | chang et al. | Deng et al. |         |
|---------|-----|--------------|--------------|--------------|-------------|---------|
| Hang Seng | 31 | 1.0974 | 1.0957 | 1.1217 | **1.0953** | 1.4181 |
| DAX 100 | 85 | 2.5424 | 2.9297 | 3.3049 | 2.5417 | **1.3190** |
| FTSE 100 | 89 | 1.1076 | 1.4623 | 1.1217 | 1.06283 | **0.8151** |
| S&P 100 | 98 | 1.9328 | 3.0696 | 3.3092 | 1.6890 | **1.4468** |
| Nikkei | 225 | 0.7961 | 0.6732 | 0.8975 | 0.6870 | **0.6179** |
| Average | - | 1.4953 | 1.8461 | 2.0483 | 1.4152 | **1.1234** |

6. **Conclusion**

In this paper we considered the portfolio selection problem under cardinality constraint which requires a predetermined number of assets to be present in the portfolio as well as bounding constraint that impose upper and lower limits on the proportions of capital invested in each asset. These real world constraints turn the problem to an NP-hard one, consequently the classical methods may not be efficient to find the optimal solution for large problem sizes.

In our work, a version of ARO is proposed to find the cardinality constrained efficient frontier. Our algorithm uses two types of mutation which modifies share indices and weights of them, separately. We also took advantage of both stochastic and chaotic variations for mutation of shares. We evaluated the performance of the proposed approach using standard data sets considered previously in the literature which are related to five major market indices containing up to 225 assets. We also compared the results with those related to some well-known heuristics proposed previously to tackle the problem. The comparison showed that our proposed ARO outperforms GA, SA and TS applied by Chang et al. (2000) to the problem, and PSO proposed by Deng et al. (2012) in most cases. Numerical results showed that by using ARO, the average error of the aforementioned test problems is reduced by approximately 20 percent of the minimum average error calculated for the above-mentioned algorithms (see Table 5).

Future work using competitive co-evolutionary genetic algorithm to deal with the problem is currently underway.

Period Portfolio Selection with Background Risk. *Entropy, 21(10), 944*. http://doi.org/10.3390/e21100944